\documentclass[11pt]{amsart}
\usepackage{amsmath,amsfonts}

\newtheorem{thm}{Theorem}[section]
\newtheorem{cor}[thm]{Corollary}
\newtheorem{lem}[thm]{Lemma}
\newtheorem{fact}[thm]{Fact}
\newtheorem{prop}[thm]{Proposition}
\newtheorem{rem}[thm]{Remark}

\newtheorem{defn}[thm]{Definition}
\newtheorem{notation}[thm]{Notation}
\numberwithin{equation}{section}

\begin{document}

\title[Smooth approximation by functions with no critical points]
{Uniform approximation of continuous functions by smooth
functions with no critical points on Hilbert manifolds}
\author{Daniel Azagra and Manuel Cepedello Boiso}
\date{December 2, 2001}

\begin{abstract}
We prove that every continuous function on a separable
infinite-dimensional Hilbert space $X$ can be uniformly
approximated by $C^\infty$ smooth functions {\em with no critical
points}. This kind of result can be regarded as a sort of very
strong approximate version of the Morse-Sard theorem. Some
consequences of the main theorem are as follows. Every two
disjoint closed subsets of $X$ can be separated by a
one-codimensional smooth manifold which is a level set of a
smooth function with no critical points; this fact may be viewed
as a nonlinear analogue of the geometrical version of the
Hahn-Banach theorem. In particular, every closed set in $X$ can
be uniformly approximated by open sets whose boundaries are
$C^\infty$ smooth one-codimensional submanifolds of $X$. Finally,
since every Hilbert manifold is diffeomorphic to an open subset
of the Hilbert space, all of these results still hold if one
replaces the Hilbert space $X$ with any smooth manifold $M$
modelled on $X$.
\end{abstract}

\maketitle

\section[Introduction]{Introduction and main results}

A fundamental result in differential topology and analysis is the
Morse-Sard theorem \cite{Sard1, Sard2}, which states that if
$f:\mathbb{R}^{n}\longrightarrow \mathbb{R}^{m}$ is a $C^r$ smooth
function, with $r>\max\{n-m, 0\}$, and $C_{f}$ stands for the set
of critical points of $f$ (that is, the points of $X$ at which the
differential of $f$ is not surjective), then the set of critical
values, $f(C_{f})$, is of (Lebesgue) measure zero in
$\mathbb{R}^{m}$. This result also holds true for smooth functions
$f:X\longrightarrow Y$ between two smooth manifolds of dimensions
$n$ and $m$ respectively.

Several authors have dealt with the question as to what extent
one can obtain a similar result for infinite-dimensional spaces
or manifolds modelled on such spaces. Let us recall some of their
results.

Smale \cite{Smale} proved that if $X$ and $Y$ are separable
connected smooth manifolds modelled on Banach spaces and
$f:X\longrightarrow Y$ is a $C^r$ Fredholm map (that is, every
differential $df(x)$ is a Fredholm operator between the
corresponding tangent spaces) then $f(C_{f})$ is meager, and in
particular $f(C_{f})$ has no interior points, provided that
$r>\max\{\textrm{index}(df(x)), 0\}$ for all $x\in X$; here
index($df(x)$) stands for the index of the Fredholm operator
$df(x)$, that is, the difference between the dimension of the
kernel of $df(x)$ and the codimension of the image of $df(x)$,
which are both finite. However, these assumptions are quite
restrictive: for instance, if $X$ is infinite-dimensional then
there is no Fredholm map $f:X\longrightarrow\mathbb{R}$. In
general, the existence of a Fredholm map $f$ from a manifold $X$
into another manifold $Y$ implies that $Y$ is
infinite-dimensional whenever $X$ is.

On the other hand, one cannot dream of extending the Morse-Sard
theorem to infinite dimensions without imposing strong
restrictions. Indeed, as shown by Kupka's counterexample
\cite{Kupka}, there are $C^\infty$ smooth functions
$f:X\longrightarrow\mathbb{R}$, where $X$ is a Hilbert space, so
that their sets of critical values $f(C_{f})$ contain intervals
and in particular have non-empty interior.

More recently, S. M. Bates has carried out a deep study
concerning the sharpness of the hypothesis of the Morse-Sard
theorem and the geometry of the sets of critical values of smooth
functions. In particular he has shown that the above $C^r$
smoothness hypothesis in the statement of the Morse-Sard theorem
can be weakened to $C^{r-1,1}$. See \cite{Bates1, Bates2, Bates3,
Bates4, Bates5}.

Nevertheless, for many applications of the Morse-Sard theorem, it
is often enough to know that any given function can be uniformly
approximated by a map whose set of critical values has empty
interior. In this direction, Eells and McAlpin established the
following theorem \cite{EellsMcAlpin}: if $X$ is a separable
Hilbert space, then every continuous function from $X$ into
$\mathbb{R}$ can be uniformly approximated by a smooth function
$f$ whose set of critical values $f(C_{f})$ is of measure zero.
This allowed them to deduce a version of this theorem for
mappings between smooth manifolds $M$ and $N$ modelled on $X$ and
a Banach space $F$ respectively, which they called an {\em
approximate Morse-Sard theorem}: every continuous mapping from
$M$ into $N$ can be uniformly approximated by a smooth function
$f:X\longrightarrow Y$ so that $f(C_{f})$ has empty interior.
However, this seemingly much more general version of the result
is a bit tricky: indeed, as they already observed
(\cite{EellsMcAlpin}, Remark 3A), when $F$ is
infinite-dimensional, the function $f$ they obtain satisfies that
$C_{f}=X$, although $f(X)$ has empty interior in $Y$.
Unfortunately, even though all the results of that paper seem to
be true, some of the proofs are not correct.

In this paper we will prove a much stronger result: if $M$ is a
$C^\infty$ smooth manifold modelled on a separable
infinite-dimensional Hilbert space $X$ (in the sequel such a
manifold will be called a Hilbert manifold), then every continuous
function on $M$ can be uniformly approximated by $C^\infty$
smooth functions {\em with no critical points}. This kind of
result might be regarded as the strongest possible one of any
class of approximate Morse-Sard theorems, when the target space is
$\mathbb{R}$.

As a by-product we also obtain the following: for every open set
$U$ in a separable Hilbert space $X$ there is a $C^{\infty}$
smooth function $f$ whose support is the closure of $U$ and so
that $f'(x)\neq 0$ for every $x\in U$. This result could be
summed up by saying that for every open subset $U$ of $X$ there
is a function $f$ whose open support is $U$ and which does not
satisfy Rolle's theorem; one should compare this result with the
main theorem from \cite{AJ} (see also the references therein).

Either of these results has in turn interesting consequences
related to smooth approximation and separation of closed sets. For
instance, every closed set in a separable Hilbert manifold $M$
can be uniformly approximated by open sets whose boundaries are
smooth one-codimensional submanifolds of $M$. Moreover, every two
disjoint closed subsets in $M$ can be separated by a smooth
one-codimensional submanifold of $M$ which is a level set of a
smooth function with no critical points. The latter may in turn be
regarded as a nonlinear analogue of the geometrical version of the
Hahn-Banach theorem.

Let us now formally state our main results.

\begin{thm}\label{main theorem}
Let $U$ be an open subset of a separable infinite-dimensional
Hilbert space $X$. Then, for all continuous functions
$f:U\longrightarrow\mathbb{R}$ and $\varepsilon:U\longrightarrow
(0,+\infty)$, there are $C^\infty$ smooth functions $\psi$ on $U$
such that $|f(x)-\psi(x)|\leq\varepsilon(x)$ and $\psi'(x)\neq 0$
whenever $x\in X$.
\end{thm}

We will prove this result in the following section. Let us now
establish the announced consequences of Theorem \ref{main
theorem}.

One could probably adapt the ideas in our proof to extend Theorem
\ref{main theorem} to the setting of Hilbert manifolds but, for
simplicity, we will instead use another approach. Indeed, bearing
in mind a fundamental result on Hilbert manifolds due to Eells
and Elworthy \cite{EE} that every separable Hilbert manifold can
be $C^\infty$ embedded as an open subset of the Hilbert space, it
is a triviality to observe that Theorem \ref{main theorem} still
holds if we replace $U$ with a a separable Hilbert manifold.

\begin{thm}\label{main theorem for manifolds}
Let $M$ be a separable Hilbert manifold. Then, for all continuous
functions $f:M\longrightarrow\mathbb{R}$ and
$\varepsilon:M\longrightarrow (0,+\infty)$, there are $C^\infty$
smooth functions $\psi:M\longrightarrow\mathbb{R}$ so that
$|f(x)-\psi(x)|\leq\varepsilon(x)$, and $d\psi(x)\neq 0$, for all
$x\in X$.
\end{thm}
\begin{proof}
According to the main theorem of \cite{EE}, there is a $C^\infty$
embedding of $M$ onto an open subset of the Hilbert space $X$.
Therefore $M$ is $C^\infty$ diffeomorphic to an open subset $U$
of $X$; let $h:U\longrightarrow M$ be such a $C^\infty$
diffeomorphism. Consider the continuous functions $g=f\circ
h:U\longrightarrow\mathbb{R}$ and $\delta=\varepsilon\circ
h:U\longrightarrow (0,+\infty)$. By Theorem \ref{main theorem}
there is a $C^\infty$ smooth function
$\varphi:U\longrightarrow\mathbb{R}$ so that $\varphi$ has no
critical points, and
    $$
    |g(y)-\varphi(y)|\leq\delta(y)
    $$
for all $y\in U$. Now define $\psi=\varphi\circ
h^{-1}:M\longrightarrow\mathbb{R}$. Since $h$ is a diffeomorphism
it is clear that $h$ takes the critical set of $\psi$ onto the
critical set of $\varphi=\psi\circ h$. But, as the latter is
empty, so is the former; that is, $\psi$ has no critical points
either. On the other hand, it is clear that
    $$
    |f(x)-\psi(x)|=|g(h^{-1}(x))-\varphi(h^{-1}(x))|\leq
    \delta(h^{-1}(x))=\varepsilon(x)
    $$
for all $x\in M$.
\end{proof}

As an easy corollary we can deduce our promised nonlinear version
of the geometrical Hahn-Banach theorem.

We will say that an open subset $U$ of a Hilbert manifold $M$ is
{\em smooth} provided that its boundary $\partial U$ is a smooth
one-codimensional submanifold of $M$.

\begin{cor}
Let $M$ be a separable Hilbert manifold. Then, for every two
disjoint closed subsets $C_1$, $C_2$ of $M$, there exists a
$C^\infty$ smooth function $\varphi:X\longrightarrow\mathbb{R}$
with no critical points, such that the level set
$N=\varphi^{-1}(0)$ is a $1$-codimensional $C^\infty$ smooth
submanifold of $M$ that separates $C_{1}$ and $C_{2}$, in the
following sense. Define $U_{1}=\{x\in M : \varphi(x)<0\}$ and
$U_{2}=\{x\in M : \varphi(x)>0\}$; then $U_1$ and $U_2$ are
disjoint $C^\infty$ smooth open sets of $M$ so that $C_{i}\subset
U_{i}$ for $i=1, 2$, and $\partial U_{1}=\partial U_{2}=N$.
\end{cor}
\begin{proof}
By Urysohn's lemma there exists a continuous function
$f:M\longrightarrow [0,1]$ so that $C_{1}\subset f^{-1}(0)$ and
$C_{2}\subset f^{-1}(1)$. Taking $\varepsilon=1/3$ and applying
Theorem \ref{main theorem for manifolds} we get a $C^\infty$
smooth function $\psi:M\longrightarrow\mathbb{R}$ which has no
critical points and is so that
    $$
    |f(x)-\psi(x)|\leq 1/3
    $$
for all $x\in M$; in particular
    $$
    C_{1}\subseteq f^{-1}(0)\subseteq \psi^{-1}(-\infty,1/2):=U_{1},
    $$
and
    $$
    C_{2}\subseteq f^{-1}(1)\subseteq \psi^{-1}(1/2, +\infty):=U_{2}.
    $$
The open sets $U_1$ and $U_2$ are smooth because their common
boundary $N=\psi^{-1}(1/2)$ is a smooth one-codimensional
submanifold of $M$ (thanks to the implicit function theorem and
the fact that $d\psi(x)\neq 0$ for all $x\in N$). In order to
obtain the result in the above form it is enough to set
$\varphi=\psi-1/2$.
\end{proof}

A trivial consequence of this result is that every closed subset
of $X$ can be uniformly approximated by smooth open subsets of
$X$. In fact,

\begin{cor}
Every closed subset of a separable Hilbert manifold $M$ can be
approximated by smooth open subsets of $M$, in the following
sense: for every closed set $C\subset M$ and every open set $W$
containing $C$ there is a $C^\infty$ smooth open set $U$ so that
$C\subset U\subseteq W$.
\end{cor}

Finally, the following result, which also implies the above
corollary, tells us that for every open set $U$ in $X$ there
always exists a function whose open support is $U$ and which does
not satisfy Rolle's theorem.

\begin{thm}\label{second main theorem}
For every open subset $U$ of a Hilbert manifold $M$ there is a
continuous function $f$ on $M$ whose support is the closure of
$U$, so that $f$ is $C^\infty$ smooth on $U$ and yet $f$ has no
critical point in $U$.
\end{thm}
\begin{proof}
For the same reasons as in the proof of Theorem \ref{main theorem
for manifolds} we may assume that $U$ is an open subset of the
Hilbert space $X=\ell_2$. Let $\varepsilon:X\longrightarrow
[0,+\infty)$ be the distance function to $X\setminus U$, that is,
    $$
    \varepsilon(x)=\textrm{dist}(x, X\setminus U)=\inf\{\|x-y\|: y\in X\setminus U\}.
    $$
The function $\varepsilon$ is continuous on $X$ and satisfies that
$\varepsilon(x)>0$ if and only if $x\in U$. According to Theorem
\ref{main theorem}, and setting $f(x)=2\varepsilon(x)$, there
exists a $C^\infty$ smooth function
$\psi:U\longrightarrow\mathbb{R}$ which has no critical points on
$U$, and such that $\varepsilon$-approximates $f$ on $U$, that is,
    $$
    |2\varepsilon(x)-\psi(x)|\leq\varepsilon(x)
    $$
for all $x\in U$. This inequality implies that
    $$
    \lim_{x\to z}\psi(x)=0
    $$
for every $z\in\partial U$. Therefore, if we set $\psi=0$ on
$X\setminus U$, the extended function
$\psi:X\longrightarrow[0,+\infty)$ is continuous on the whole of
$X$, is $C^\infty$ smooth on $U$ and has no critical points on
$U$. On the other hand, $\psi(x)\geq\varepsilon(x)>0$ for all
$x\in U$, hence the support of $\psi$ is $\overline{U}$.
\end{proof}

\section[Proofs]{Proof of the main result}

The main ideas behind the proof of Theorem \ref{main theorem} are
as follows. First we use a perturbed smooth partition of unity to
approximate the given continuous function $f$. The summands of
this perturbed partition of unity are functions supported on
scalloped balls and carefully constructed in such a way that the
critical points of the approximating sum $\varphi$ are kept under
control. More precisely, those critical points consist of a
sequence of compact sets $K_n$ that are suitably isolated in
pairwise disjoint open sets $U_n$ of small diameter so that the
oscillation of both $f$ and $\varphi$ on $U_n$ is small as well.

Then we have to eliminate all of those critical points without
losing much of the approximation. To this end we compose the
approximating function $\varphi$ with a sequence of deleting
diffeomorphisms $h_{n}:X\longrightarrow X\setminus K_{n}$ which
extract each of the compact sets of critical points $K_n$ and
restrict to the identity outside each of the open sets $U_{n}$.
The infinite composition of deleting diffeomorphisms with our
function, $\psi=\varphi\circ \bigcirc_{n=1}^{\infty} h_{n}$, is
locally finite, in the sense that only a finite number of
diffeomorphisms are acting on some neighborhood of each point,
while all the rest restrict to the identity on that neighborhood.
In this way we obtain a smooth function $\psi$ which has no
critical points, and which happens to approximate the function
$\varphi$ (which in turn approximates the original $f$) because
the perturbation brought on $\varphi$ by that infinite
composition is not very important: indeed, recall that each $h_n$
restricts to the identity outside the set $U_n$ (on which
$\varphi$ has a small oscillation), and the $U_n$ are pairwise
disjoint.

We will make the proof of Theorem \ref{main theorem} in the case
of a constant $\varepsilon>0$ so as to avoid bearing an
unnecessary burden of notation. Later on we will briefly explain
what additional technical precautions must be taken in order to
deduce the general form of this result (see Remark \ref{remark for
e(x)}).

The following proposition shows the existence of a function
$\varphi$ with the above properties. Recall that $C_{\varphi}$
stands for the set of critical points of $\varphi$.

\begin{prop}\label{existence of varphi}
Let $U$ be an open subset of the separable Hilbert space $X$. Let
$f:U\longrightarrow\mathbb{R}$ be a continuous function on $X$,
and $\varepsilon>0$. Then there exist a $C^\infty$ smooth function
$\varphi:U\longrightarrow\mathbb{R}$, a sequence $(K_{n})$ of
compact sets, a sequence $(U_{n})$ of open sets, and a sequence
$(B(y_{n}, r_{n}))$ of open balls which are contained in $U$ and
whose union covers $U$, such that:
\begin{itemize}
\item[{(a)}] $C_{\varphi}\subseteq\bigcup_{n=1}^{\infty} K_{n}$;
\item[{(b)}] $K_{n}\subset U_{n}\subseteq B(y_{n}, r_{n})$ for all
$n\in\mathbb{N}$, and $U_{n}\cap U_{m}=\emptyset$ whenever $n\neq
m$;
\item[{(c)}] $|\varphi(x)-f(y)|\leq 2\varepsilon$ for all
$x,y\in B(y_{n},r_{n})$,  $n\in\mathbb{N}$;
\item[{(d)}] for every $x\in U$ there exist an open neighborhood
$V_{x}$ of $x$ and some $n_{x}\in\mathbb{N}$ such that $V_{x}\cap
U_{m}=\emptyset$ for all $m>n_{x}$.
\end{itemize}
\end{prop}

The following theorem ensures the existence of the diffeomorphisms
$h_{n}$.

\begin{thm}\label{removing compact sets}
Let $X$ be an infinite-dimensional Hilbert space. Then, for every
compact set $K$ and every open subset $U$ of $X$ with $K\subset
U$, there exists a $C^\infty$ smooth diffeomorphism
$h:X\longrightarrow X\setminus K$ so that $h$ restricts to the
identity outside $U$.
\end{thm}
\noindent This result may be regarded, in the Hilbert case, as a
(rather technical, but crucial to our purposes) improvement of
some known results on smooth negligibility of compact sets (see
\cite{ADo, Do}; there $h$ is known to be the identity only
outside a ball containing $K$).

Assume for a while that Proposition \ref{existence of varphi} and
Theorem \ref{removing compact sets} are already established, and
let us see how we can deduce Theorem \ref{main theorem}.

\medskip

\begin{center}
{\bf Proof of Theorem \ref{main theorem}}
\end{center}

For a given continuous function $f$ and a number $\varepsilon>0$,
take a function $\varphi$ and sequences $(K_{n})$ and $(U_{n})$
with the properties of Proposition \ref{existence of varphi}. For
each compact set $K_{n}$ and each open set $U_{n}$, use Theorem
\ref{removing compact sets} to find a $C^{\infty}$ diffeomorphism
$h_{n}:X\longrightarrow X\setminus K_{n}$ so that $h_{n}(x)=x$ if
$x\notin U_{n}$. Note that, since the $U_j$ contain the $K_j$ and
are pairwise disjoint,
    $$
    h_{n}(x)\notin\bigcup_{j=1}^{\infty}K_{j} \supseteq C_{\varphi} \eqno (1)
    $$
for all $x\in U$, $n\in\mathbb{N}$. Define then
$\psi:U\longrightarrow\mathbb{R}$ by
    $$
    \psi=\varphi\circ\bigcirc_{n=1}^{\infty} h_{n}.
    $$
This formula makes sense and the function $\psi$ is $C^\infty$
smooth because the infinite composition is in fact locally finite.
Indeed, for a given $x\in U$, according to Proposition
\ref{existence of varphi}(d), we can find an open neighborhood
$V_{x}$ of $x$ and some $n_{x}\in\mathbb{N}$ so that $V_{x}\cap
U_{m}=\emptyset$ for all $m>n_{x}$; hence $h_{m}(y)=y$ for all
$y\in V_{x}$ and $m>n_{x}$, and therefore
    $$
    \psi(y)=\varphi\circ h_{n_{x}}\circ h_{n_{x}-1}\circ ...\circ
    h_{2}\circ h_{1}(y) \eqno(2)
    $$
for all $y\in V_{x}$. The derivative $\psi'(y)$ is given by
    $$
    \psi'(y)=\varphi'\bigl(\bigcirc_{j=1}^{n_{x}}h_{j}(y)\bigr)\circ Dh_{n_{x}}
    \bigl(\bigcirc_{j=1}^{n_{x}-1}h_{j}(y)\bigr)\circ ...\circ
    Dh_{2}(h_{1}(y))\circ Dh_{1}(y) \eqno(3)
    $$
for all $y\in V_{x}$. Since $U_{n}\subseteq X\setminus U_{m}$ for
$n\neq m$, we have that $h_{m}$ is the identity on $U_{n}$, and
therefore $Dh_{m}(x)=I$ (the identity isomorphism of $\ell_{2}$)
for all $x\in U_{n}$. By the continuity of $Dh_{n}$ it follows
that $Dh_{m}(x)=I$ for all $x\in\overline{U}_{n}$, if $m\neq n$.
This implies that, for $y\in\overline{U_{n}}\cap V_{x}$, all the
differentials $Dh_{j}(z)$ in $(3)$ are the identity, except
perhaps for $j=n$. Hence we get that either
    $$
    \psi'(y)=\varphi'(h_{n}(y))\circ Dh_{n}(y), \, \textrm{ and }\,
    \psi(y)=\varphi(h_{n}(y)), \eqno(4)
    $$
if $y$ belongs to some $\overline{U}_{n}$, or else
    $$
    \psi'(y)=\varphi'(y),  \, \textrm{ and }\,
    \psi(y)=\varphi(y), \eqno(5)
    $$
when $y\notin\bigcup_{n=1}^{\infty}\overline{U}_{n}$.

Now we can easily check that $C_{\psi}=\emptyset$. Take $x\in U$.
If we are in the case that $x\in\overline{U}_{n}$ for some $n$
then
    $
    \psi'(x)=\varphi'(h_{n}(x))\circ Dh_{n}(x)\neq 0,
    $
because $D h_{n}(x)$ is a linear isomorphism and, according to
$(1)$ above, $\varphi'(h_{n}(x))\neq 0$. Otherwise we have that
$x\notin\bigcup_{n=1}^{\infty}\overline{U}_{n}\supseteq
C_{\varphi}$, so $\psi'(x)=\varphi'(x)\neq 0$ trivially.

It only remains to check that $\psi$ still approximates $f$. As
before, for a given $x\in U$, either $\psi(x)=\varphi(x)$ or
$\psi(x)=\varphi(h_{n}(x))$ for some $n$ (with $x\in U_{n}$). In
the first case, from Proposition \ref{existence of varphi}(c) we
get that $|\psi(x)-f(x)|\leq 2\varepsilon$. In the second case,
bearing in mind that $h_{n}(x)\in U_{n}\subseteq B(y_{n},r_{n})$,
and for the same reason, we have that
    $$
    |\psi(x)-f(x)|=|\varphi(h_{n}(x))-f(x)|\leq 2\varepsilon;
    $$
in either case we obtain that $|\psi(x)-f(x)|\leq 2\varepsilon$.


\medskip

\begin{center}
{\bf Proof of Proposition \ref{existence of varphi}}
\end{center}

We will assume that $U=X$, since the proof is completely analogous
in the case of a general open set. One only has to take some (easy
but rather rambling) technical precautions in order to make sure
that the different balls considered in the argument are in $U$.

Let $B(x, r)$ and $\overline{B}(x,r)$ stand for the open ball and
closed ball, respectively, of center $x$ and radius $r$, with
respect to the usual hilbertian norm $\|\cdot\|$ of $X$.

Let $f:X\longrightarrow\mathbb{R}$ be a continuous function, and
let $\varepsilon>0$. By continuity, for every $x\in X$ there
exists $\delta_{x}>0$ so that $|f(y)-f(x)|\leq \varepsilon/4$
whenever $y\in B(x, 2\delta_{x})$. Since $X=\bigcup_{x\in X}
B(x,\delta_{x}/2)$ is separable, there exists a countable
subcovering,
    $$
    X=\bigcup_{n=1}^{\infty}B(x_{n}, s_{n}/2),
    $$
where $s_{n}=\delta_{x_{n}}$, for some sequence of centers
$(x_{n})$. By induction we can choose a sequence of {\em linearly
independent} vectors $(y_{n})$, with $y_{n}\in B(x_{n}, s_{n}/2)$,
so that
    $$
    X=\bigcup_{n=1}^{\infty}B(y_{n}, s_{n}).
    $$
Moreover, we have that
    $$
    |f(y)-f(y_{n})|\leq\varepsilon/2
    $$
provided $\|y-y_{n}\|\leq \frac{3}{2} s_{n}$, as is immediately
checked.
\begin{notation}
{\em In the sequel $\mathcal{A}[z_{1},...,z_{k}]$ stands for the
affine subspace spanned by a finite sequence of points $z_{1},
..., z_{k}\in X$.}
\end{notation}
The following lemma shows that we can slightly move the radii
$s_{n}$ so that, for any finite selection of centers $y_{n}$, the
spheres that are the boundaries of the balls $B(y_{n}, s_{n})$
have empty intersection with the affine subspace spanned by those
centers.
\begin{lem}\label{controlled intersection of spheres}
We can find a sequence of positive numbers $(r_{n})$ with
$s_{n}\leq r_{n}\leq \frac{3}{2} s_{n}$ so that, if we denote
$S_{n}=\partial B(y_{n}, r_{n})$ then,
\begin{itemize}
\item[{(i)}] for each finite sequence of positive integers
$k_{1}<k_{2}<...<k_{m}$,
    $$
    \mathcal{A}[y_{k_{1}}, ..., y_{k_{m}}]\cap S_{k_{1}}\cap ...\cap
    S_{k_{m}}=\emptyset.
    $$
\item[{(ii)}] for any $n,k\in\mathbb{N}$, $y_{n}\notin S_{k}$.
\end{itemize}
\end{lem}
\begin{proof}
We will define the $r_{n}$ inductively.

For $n=1$ we may take $r_{1}\in [s_{1}, \frac{3}{2}s_{1}]$ so
that $r_{1}$ does not belong to the countable set
$\{\|y_{1}-y_{k}\| : k\in\mathbb{N}\}$; this means that
$y_{k}\notin S_{1}$ for any $k\in\mathbb{N}$. On the other hand,
it is obvious that $\{y_{1}\}\cap S_{1}=\emptyset$.

Assume now that $r_{1}, ..., r_{n}$ have already been chosen in
such a way that the spheres $S_{1}$, ..., $S_{n}$ satisfy $(i)$
and $(ii)$, and let us see how we can find $r_{n+1}$. For any
finite sequence of integers $0<k_{1}< ... <k_{j}\leq n+1$, let us
denote
    $$
    \mathcal{A}_{k_{1},...,k_{j}}=\mathcal{A}[y_{k_{1}}, ...,
    y_{k_{j}}].
    $$
For simplicity, and up to a suitable translation (which obviously
does not affect our problem), we may assume that $y_{n+1}=0$, so
that $\mathcal{A}_{k_{1},...,k_{m},n+1}$ is the $m$-dimensional
vector subspace of $X$ spanned by $y_{k_{1}}$, ..., $y_{k_{m}}$.
Now, for each finite sequence of integers $0<k_{1}<...<k_{m}\leq
n$, consider the map $F_{k_{1}, ...,
k_{m}}:\mathcal{A}_{k_{1},...,k_{m}, n+1}\longrightarrow
\mathbb{R}^{m}$ defined by
    $$
    F_{k_{1}, ..., k_{m}}(x)=\big(\|x-y_{k_{1}}\|^{2}- {r_{k_{1}}}^{2}, ...,
    \|x-y_{k_{m}}\|^{2}- {r_{k_{m}}}^{2}\big).
    $$
Note that
    $$
    D F_{k_{1}, ..., k_{m}}(x)=\big(2(x-y_{k_{1}}), ..., 2(x-y_{k_{m}})\big)
    $$
and therefore $\textrm{rank}\big(D F_{k_{1}, ...,
k_{m}}(x)\big)<m$ if and only if $x\in
\mathcal{A}_{k_{1},...,k_{m}}$. By the induction assumption we
know that
    $$
    S_{k_{1}}\cap ... \cap S_{k_{m}}\cap \mathcal{A}_{k_{1},...,k_{m}}=\emptyset,
    $$
hence it is clear that $\textrm{rank}\big(D F_{k_{1}, ...,
k_{m}}(x)\big)=m$ for all $x\in S_{k_{1}}\cap ... \cap
S_{k_{m}}\cap \mathcal{A}_{k_{1},...,k_{m}, n+1}$. This implies
that
    $$
    M_{k_{1}, ..., k_{m}}:=S_{k_{1}}\cap ... \cap
S_{k_{m}}\cap \mathcal{A}_{k_{1},...,k_{m}, n+1}
    $$
is a compact $m-m=0$-dimensional submanifold of
$\mathcal{A}_{k_{1},...,k_{m}, n+1}$, and in particular
$M_{k_{1}, ..., k_{m}}$ consists of a finite number of points (in
fact two points, but we do not need to know this). Therefore
    $$
    M=\bigcup M_{k_{1}, ..., k_{m}}
    $$
(where the union is taken over all the finite sequences of
integers $0<k_{1}<...<k_{n}\leq n$) is a finite set as well. Now
we have that
    $$I:=\big[s_{n+1}, \frac{3}{2}s_{n+1}\big]\setminus
    \Big(\{\|z\| : z\in M\}\cup\{\|y_{j}\| :
    j\in\mathbb{N}\}\Big)
    $$
is an uncountable subset of the real line, so we can find a number
$r_{n+1}\in I$. With this choice it is clear that
    $$
    S_{k_{1}}\cap ... \cap
S_{k_{m}}\cap S_{n+1}\cap \mathcal{A}_{k_{1},...,k_{m}, n+1}=
M_{k_{1},...,k_{m}}\cap S_{n+1}=\emptyset
    $$
for all finite sequences of integers $0<k_{1}<...<k_{m}<n+1$, and
also $y_{j}\notin S_{n+1}=\partial B(0, r_{n+1})$ for all
$j\in\mathbb{N}$. Therefore the spheres $S_{1}$, ..., $S_{n}$,
$S_{n+1}$ satisfy $(i)$ and $(ii)$ as well. By induction the
sequence $(r_{n})$ is thus well defined.
\end{proof}

\medskip
Since $s_{n}\leq r_{n}\leq \frac{3}{2} s_{n}$ for all $n$, it is
clear that the new balls $B(y_{n},r_{n})$ keep the two important
properties of the old balls $B(y_{n}, s_{n})$, namely,
    $$
    X=\bigcup_{n=1}^{\infty}B(y_{n}, r_{n}), \eqno (6)
    $$
and
    $$
    |f(y)-f(y_{n})|\leq\varepsilon/2 \text{ whenever $\|y-y_{n}\|\leq r_{n}$.} \eqno (7)
    $$

Now we define the scalloped balls $B_{n}$ that are the basis for
our perturbed partition of unity: set $B_{1}=B(y_{1},r_{1})$, and
for $n\geq 2$ define
    $$
    B_{n}=B(y_{n}, r_{n})\setminus \Bigl(\bigcup_{j=1}^{n-1}
    \overline{B}(y_{j}, \lambda_{n} r_{j})\Bigr);
    $$
where $1/2<\lambda_{2}<\lambda_{3}< ...
<\lambda_{n}<\lambda_{n+1}< ... <1$, with $\lim_{n\to\infty}\lambda_{n}=1$.
The $\lambda_{n}$ are to be fixed later on.

Taking into account that $\lim_{n\to\infty}\lambda_{n}=1$, it is
easily checked that the $B_{n}$ form a locally finite open
covering of $X$, with the nice property that
    $$
    |f(y)-f(y_{n})|\leq\varepsilon/2
    \text{ whenever $y\in B_{n}$.}
    $$

Next, pick a $C^{\infty}$ smooth function
$g_{1}:\mathbb{R}\longrightarrow [0,1]$ so that:
\begin{itemize}
\item[{(i)}] $g_{1}(t)=1$ for $t\leq 0$,
\item[{(ii)}] $g_{1}(t)=0$ for $t\geq {r_{1}}^{2}$,
\item[{(iii)}] $g_{1}'(t)<0$ if $0<t< {r_{1}}^{2}$;
\end{itemize}
and define then $\varphi_{1}:X\longrightarrow\mathbb{R}$ by
    $$
    \varphi_{1}(x)=g_{1}(\|x-y_{1}\|^{2})
    $$
for all $x\in X$. Note that $\varphi_{1}$ is a $C^{\infty}$ smooth
function whose open support is $B_{1}$, and $B_{1}\cap
C_{\varphi_{1}}=\{y_{1}\}$, that is, $y_{1}$ is the only critical
point of $\varphi_{1}$ that lies inside $B_{1}$.

Now, for $n\geq 2$, pick $C^\infty$ smooth functions
$\theta_{(n,j)}:\mathbb{R}\longrightarrow [0,1]$, $j=1, ..., n$,
with the following properties. For $j=1, ..., n-1$,
$\theta_{(n,j)}$ satisfies that
\begin{itemize}
\item[{(i)}] $\theta_{(n,j)}(t)=0$ for $t\leq(\lambda_{n}r_{j})^{2}$,
\item[{(ii)}] $\theta_{(n,j)}(t)=1$ for $t\geq {r_{j}}^{2}$,
\item[{(iii)}] $\theta_{(n,j)}'(t)>0$ if $(\lambda_{n}r_{j})^{2}<t<{r_{j}}^{2}$;
\end{itemize}
while for $j=n$ the function $\theta_{(n,n)}$ is such that
\begin{itemize}
\item[{(i)}] $\theta_{(n,n)}(t)=1$ for $t\leq 0$,
\item[{(ii)}] $\theta_{(n,n)}(t)=0$ for $t\geq {r_{n}}^{2}$,
\item[{(iii)}] $\theta_{(n,n)}'(t)<0$ if $0<t< {r_{n}}^{2}$.
\end{itemize}
Then define the function $g_{n}:\mathbb{R}^{n}\longrightarrow
[0,1]$ as
    $$
    g_{n}(t_{1}, ..., t_{n})=\prod_{i=1}^{n}\theta_{(n,i)}(t_{i})
    $$
for all $t=(t_{1}, ..., t_{n})\in\mathbb{R}^{n}$. This function
is clearly $C^\infty$ smooth on $\mathbb{R}^{n}$ and satisfies
the following properties:
\begin{itemize}
\item[{(i)}] $g_{n}(t_{1},...,t_{n})>0$ if and only if
$t_{j}>(\lambda_{n}r_{j})^{2}$ for all $j=1, ..., n-1$, and
$t_{n}< {r_{n}}^{2}$; and $g_{n}$ vanishes elsewhere;
\item[{(ii)}] $g_{n}(t_{1},...,t_{n})=\theta_{(n,n)}(t_{n})$ whenever $t_{j}\geq
{r_{j}}^{2}$ for all $j=1, ..., n-1$;
\item[{(iii)}] $\nabla g_{n}(t_{1},...,t_{n})\neq 0$ provided
$(\lambda_{n}r_{j})^{2}<t_{j}$ for all $j=1, ..., n-1$, and
$0<t_{n}<{r_{n}}^{2}$.
\end{itemize}
Moreover, under the same conditions as in (iii) just above we
have that
    $$
    \frac{\partial g_{n}}{\partial t_{n}}(t_{1},...,t_{n})=
    \frac{\partial \theta_{(n,n)}}{\partial t_{n}}(t_{n})
    \prod_{i=1}^{n-1}\theta_{(n,i)}(t_{i})<0, \eqno(8)
    $$
since no function in this product vanishes on the specified set,
while for $j<n$, according to the corresponding properties of the
functions $\theta_{(n,j)}$ we have that
    $$
    \frac{\partial g_{n}}{\partial t_{j}}(t_{1},...,t_{n})=
    \frac{\partial \theta_{(n,j)}}{\partial t_{j}}(t_{j})
    \prod_{i=1, i\neq j}^{n}\theta_{(n,i)}(t_{i})>0. \eqno(9)
    $$
If we are not in the conditions of (iii) then the corresponding
inequalities do still hold but are not strict.

Let us now define $\varphi_{n}:X\longrightarrow [0,1]$ by
    $$
    \varphi_{n}(x)=g_{n}(\|x-y_{1}\|^{2}, ... , \|x-y_{n}\|^{2}).
    $$
It is clear that $\varphi_{n}$ is a $C^{\infty}$ smooth function
whose open support is precisely the scalloped ball $B_{n}$.

As above, let us denote by $C_{\varphi_{n}}$ the critical set of
$\varphi_{n}$, that is,
$$
C_{\varphi_{n}}=\{x\in X :
\varphi_{n}'(x)=0\}.
$$
Since our norm $\|\cdot\|$ is hilbertian we
have that, if $x\in C_{\varphi_{n}}\cap B_{n}$, then $x$ belongs
to the affine span of $y_{1}, ..., y_{n}$. Indeed, if $x\in
B_{n}$,
    $$
    \varphi_{n}'(x)=\sum_{j=1}^{n}\frac{\partial g_{n}}{\partial t_{j}}
    (\|x-y_{1}\|^{2}, ... , \|x-y_{n}\|^{2})\,2(x-y_{j})=0, \eqno(10)
    $$
which (taking into account (8) and the fact that the $y_{j}$ are
all linearly independent) means that $x$ is in the affine span of
$y_{1}, ..., y_{n}$. Here, as is usual, we identify the Hilbert
space $X$ with its dual $X^*$, and we make use of the fact that
the derivative of the function $x\mapsto \|x\|^{2}$ is the
mapping $x\mapsto 2x$.

Similarly, by using (8) it can be shown that $x\in
C_{\varphi_{1}+\cdots+\varphi_{m}}\cap (B_{1}\cup ...\cup B_{m})$
implies that $x$ belongs to the affine span of $y_{1}, ...,
y_{m}$.

\medskip

In order that our approximating function has a small critical set
we cannot use the standard approximation provided by the partition
of unity associated with the functions
$(\varphi_{j})_{i\in\mathbb{N}}$, namely
    $$
    x\mapsto \frac{\sum_{n=1}^{\infty}\alpha_{n}\varphi_{n}(x)}
    {\sum_{n=1}^{\infty}\varphi_{n}(x)},
    $$
where $\alpha_{n}=f(y_{n})$. Indeed, such a function would have a
huge set of critical points since it would be constant (equal to
$\alpha_{n}$) on a lot of large places (at least on each $B_{n}$
minus the union of the rest of the $B_j$). Instead, we will
modify this standard approximation by letting the $\alpha_{n}$ be
functions (and not mere numbers) of very small oscillation and
with only one critical point (namely $y_{n}$). So, for every
$n\in\mathbb{N}$ let us pick a $C^{\infty}$ smooth real function
$a_{n}:[0, +\infty)\longrightarrow\mathbb{R}$ with the following
properties:
\begin{itemize}
\item[{(i)}] $a_{n}(0)=f(y_{n})$;
\item[{(ii)}] $a_{n}'(t)<0$ whenever $t>0$;
\item[{(iii)}] $|a_{n}(t)-a_{n}(0)|\leq\varepsilon/2$ for all $t\geq 0$;
\end{itemize}
and define $\alpha_{n}:X\longrightarrow\mathbb{R}$ by
    $$
    \alpha_{n}(x)=a_{n}(\|x-y_{n}\|^{2})
    $$
for every $x\in X$. It is clear that $\alpha_{n}$ is a $C^\infty$
smooth function on $X$ whose only critical point is $y_{n}$.
Besides,
$$
|\alpha_{n}(x)-f(y_{n})|\leq\varepsilon/2 \text{ for all $x\in X$}.
$$

Now we can define our approximating function
$\varphi:X\longrightarrow\mathbb{R}$ by
    $$
    \varphi(x)=\frac{\sum_{n=1}^{\infty}\alpha_{n}(x)\varphi_{n}(x)}
    {\sum_{n=1}^{\infty}\varphi_{n}(x)}
    $$
for every $x\in X$. Since the sums are locally finite, it is clear
that $\varphi$ is a well-defined $C^\infty$ smooth function.

\begin{fact}\label{good approximation}
The function $\varphi$ approximates $f$ nicely. Namely, we have
that
\begin{itemize}
\item[{(i)}] $|\varphi(x)-f(x)|\leq\varepsilon$ for all $x\in X$, and
\item[{(ii)}] $|\varphi(y)-f(x)|\leq 2\varepsilon$ for all $x, y\in
B(y_{n}, r_{n})$ and each $n\in\mathbb{N}$.
\end{itemize}
\end{fact}
\begin{proof}
Indeed, for every $n$ we have that
$|\alpha_{n}(x)-f(y_{n})|\leq\varepsilon/2$ for all $x\in X$. On
the other hand, by $(7)$ above we know that $|f(x)-f(y_{n})|\leq
\varepsilon/2$ whenever $x\in B(y_{n}, r_{n})$. Then, by the
triangle inequality, it follows that
    $$
    |\alpha_{n}(x)-f(x)|\leq\varepsilon \eqno(11)
    $$
whenever $x\in B(y_{n}, r_{n})$. In the same way we deduce that
    $$
    |\alpha_{m}(x)-f(y_{n})|\leq\varepsilon \eqno(12)
    $$
whenever $x\in B(y_{n}, r_{n})\cap B(y_{m}, r_{m})$. Since
$\varphi_{m}(y)=0$ when $y\notin B(y_{m}, r_{m})$, from $(11)$ we
get that
    $$
    |\varphi(x)-f(x)|=
    \bigg|\frac{\sum_{m=1}^{\infty}(\alpha_{m}(x)-f(x))\varphi_{m}(x)}
    {\sum_{m=1}^{\infty}\varphi_{m}(x)}\bigg|\leq
    \frac{\sum_{m=1}^{\infty}\varepsilon \varphi_{m}(x)}
    {\sum_{m=1}^{\infty}\varphi_{m}(x)}=\varepsilon
    $$
for all $x\in X$, which shows $(i)$. Similarly, we deduce from
$(12)$ that
    $$
    |\varphi(y)-f(y_{n})|=
    \bigg|\frac{\sum_{m=1}^{\infty}(\alpha_{m}(y)-f(y_{n}))\varphi_{m}(y)}
    {\sum_{m=1}^{\infty}\varphi_{m}(y)}\bigg|\leq
    \frac{\sum_{m=1}^{\infty}\varepsilon \varphi_{m}(y)}
    {\sum_{m=1}^{\infty}\varphi_{m}(y)}=\varepsilon
    $$
for every $y\in B(y_{n}, r_{n})$, which, combined with the fact
that $|f(x)-f(y_{n})|\leq\varepsilon/2$ for $x\in B(y_{n},
r_{n})$, yields that
    $$
    |\varphi(y)-f(x)|\leq \varepsilon+\varepsilon/2,
    $$
for every $x, y\in B(y_{n}, r_{n})$, so $(ii)$ is satisfied as
well.
\end{proof}

\medskip

Now let us have a look at the derivative of $\varphi$. To this
end let us introduce the auxiliary functions $f_{n}$ defined by
    $$
    f_{n}(x)=\frac{\sum_{k=1}^{n}\alpha_{k}(x)\varphi_{k}(x)}
    {\sum_{k=1}^{n}\varphi_{k}(x)}, \, \textrm{ for all }\,
    x\in \bigcup_{i=1}^{n}B_{i}.
    $$
Notice that $\varphi$ can be expressed as
    $$
    \varphi(x)=\lim_{n\to\infty}f_{n}(x),
    $$
that the domains of the $f_{n}$ form an increasing tower of open
sets whose union is $X$, and that each $f_{n}$ restricts to
$f_{n-1}$ on $\bigcup_{i=1}^{n-1}B_{i}\setminus B_{n}$. Moreover,
for each $x\in X$ there is some open neighborhood $V_{x}$ of $x$
and some $n_{x}\in\mathbb{N}$ so that $\varphi(y)=f_{n_{x}}(y)$
for all $y\in V_{x}$. In fact we have that
    $$
    \varphi(x)=f_{n}(x) \hspace{0.3cm} \textrm{for all} \hspace{0.3cm}
    x\in V_{n}:=\big(\bigcup_{j=1}^{n}B_{j}\big)\setminus
    \big(\bigcup_{i=n+1}^{\infty}\overline{B}_{i}\big),
    $$
for every $n$, the $V_{n}$ are open, $V_{n}\subseteq V_{n+1}$,
and $\bigcup_{i=1}^{\infty}V_{i}=X$, because the covering of $X$
formed by the $B_j$ is locally finite.

Hence, by looking at the derivatives of the functions $f_{n}$ we
will get enough information about the derivative of $\varphi$.

If $x\in\bigcup_{j=1}^{n}B_{j}$ then the expression for the
derivative of $f_{n}$ is given by
    $$
    f_{n}'(x)=\frac{\sum_{j=1}^{n}[\alpha_{j}'(x)\varphi_{j}(x)+
    \alpha_{j}(x)\varphi_{j}'(x)]\sum_{i=1}^{n}\varphi_{i}(x)-
    \sum_{j=1}^{n}\varphi_{j}'(x)\sum_{i=1}^{n}\alpha_{i}(x)\varphi_{i}(x)}
    {(\sum_{j=1}^{n}\varphi_{j}(x))^{2}}.
    $$
Therefore, for $x\in\bigcup_{j=1}^{n}B_{j}$ we have that
$f_{n}'(x)=0$ if and only if
    $$
    \sum_{j=1}^{n}\sum_{i=1}^{n}\varphi_{i}(x)
    \Bigl[\alpha_{j}'(x)\varphi_{j}(x)+\bigl(\alpha_{j}(x)-\alpha_{i}(x)\bigr)\varphi_{j}'(x)\Bigr]
    =0. \eqno(13)
    $$
By inserting the expressions for the derivatives of $\varphi_{j}$
and $\alpha_{j}$ in equation $(13)$, we can express the condition
$f_{n}'(x)=0$ as a nontrivial linear dependence link on the
vectors $(x-y_{j})$, which yields that $x$ is in the affine span
of the points $y_{1}, ..., y_{n}$. Indeed, we are going to prove
the following.
\begin{fact}\label{the critical set is locally finite
dimensional} If $x\in C_{f_{n}}\cap B_{n}$, then $x\in
\mathcal{A}_{n}:=\mathcal{A}[y_{1}, ..., y_{n}]$. Moreover, for
each $n\in\mathbb{N}$ and for every finite sequence of positive
integers $k_{1}<k_{2}<...<k_{m}<n$ we have that
    $$
    C_{f_{n}}\cap\big(B_{n}\setminus\bigcup_{j=1}^{m}B_{k_{j}}\big)
    \subseteq\mathcal{A}\big[\{y_{1}, ..., y_{n}\}
    \setminus\{y_{k_{1}}, ..., y_{k_{m}}\}\big].
    $$
\end{fact}
\begin{proof}
As above, in all the subsequent calculations, we will identify the
Hilbert space $X$ with its dual $X^*$, and the derivative of
$\|\cdot\|^{2}$ with the mapping $x\mapsto 2x$. To save notation,
let us simply write
    $$
    \frac{\partial g_{n}}{\partial
    t_{j}}(\|x-y_{1}\|^{2},...,\|x-y_{n}\|^{2})=\mu_{(n,j)},
    $$
and
    $$
    a_{j}'(\|x-y_{j}\|^{2})=\eta_{j}.
    $$

Notice that, according to $(8)$ and $(9)$ above, $\mu_{(n,j)}\geq
0$ for $j=1, ..., n-1$, while $\mu_{(n,n)}\leq 0$; and
$\mu_{(n,n)}\neq 0$ provided $x\in B_{n}$ and $x\neq y_{n}$; on
the other hand it is clear that $\eta_{j}<0$ for all $j$ unless
$x=y_{j}$ (in which case $\eta_{j}=0$).

Assuming $x\in C_{f_{n}}\cap B_{n}$, and taking into account the
expression $(10)$ for $\varphi_{j}'(x)$ and the fact that
$\alpha_{j}'(x)=2\eta_{j}(x-y_{j})$, we can write condition $(13)$
above in the form
    $$
    2\sum_{j=1}^{n}\sum_{i=1}^{n}\varphi_{i}(x)
    \Bigl[\eta_{j}\varphi_{j}(x)\,(x-y_{j})+\bigl(\alpha_{j}(x)-\alpha_{i}(x)\bigr)
    \sum_{\ell=1}^{j}\mu_{(j,\ell)}\,(x-y_{\ell})\Bigr]
    =0,
    $$
which in turn is equivalent (taking the common factors of each
$(x-y_{j})$ together) to the following one
    $$
    \sum_{j=1}^{n}\Biggl[
    \eta_{j}\varphi_{j}(x)\sum_{i=1}^{n}\varphi_{i}(x) +
    \sum_{k=j}^{n}\Bigl(\sum_{i=1}^{n}\bigl(\alpha_{k}(x)-
    \alpha_{i}(x)\bigr)\varphi_{i}(x)\Bigr)
    \mu_{(k,j)}\Biggr]\,(x-y_{j})=0. \eqno(14)
    $$
Now notice that, if we can prove that at least one of the
expressions multiplying the $(x-y_{j})$ does not vanish then we
are done; indeed, we will have that the vectors $x-y_{1}$, ...,
$x-y_{n}$ are linearly dependent, which means that $x$ belongs to
the affine span of the points $y_{1}, ..., y_{n}$.

So let us check that not all of those expressions in $(14)$
vanish. In fact we are going to see that at least one of the
terms is strictly negative. We can obviously assume that $x$ is
not any of the points $y_{1}, ..., y_{n}$ (which are already in
$\mathcal{A}_{n}$). In this case we have that $\mu_{(n,n)}<0$ and
$\eta_{j}<0$ for all $j=1, ..., n$. For simplicity, we will only
make the argument in the case $n=3$; giving a proof in a more
general case would be as little instructive as tedious to read.

Let us first assume that $\varphi_{j}(x)\neq 0$ for $j=1,2,3$. We
begin by looking at the term that multiplies $(x-y_{3})$ in
$(14)$, that is
    $$
    \beta_{3}:=\eta_{3}\varphi_{3}(x)\sum_{i=1}^{3}\varphi_{i}(x) +
    \sum_{i=1}^{3}\bigl(\alpha_{3}(x)-
    \alpha_{i}(x)\bigr)\varphi_{i}(x)
    \mu_{(3,3)}.
    $$
If $\sum_{i=1}^{3}\bigl(\alpha_{3}(x)-\alpha_{i}(x)\bigl) \varphi_{i}(x)\geq
0$ we are done, since in this case we easily see that
$\beta_{3}<0$ (remember that $\mu_{(3,3)}\leq 0$, $\eta_{3}<0$,
and $\varphi_{3}(x)>0$). Otherwise we have that
$$
\sum_{i=1}^{3}\Bigl(\alpha_{3}(x)-\alpha_{i}(x)\Bigr) \varphi_{i}(x)<0,
$$
and then we look at the term $\beta_{2}$ multiplying $(x-y_{2})$
in $(14)$, namely,
    $$
    \beta_{2}:=\eta_{2}\varphi_{2}(x)\sum_{i=1}^{3}\varphi_{i}(x) +
    \sum_{k=2}^{3}\Biggl(\sum_{i=1}^{3}\Bigl(\alpha_{k}(x)-
    \alpha_{i}(x)\Bigr)\varphi_{i}(x)\Biggr)
    \mu_{(k,2)}.
    $$
Now, since $\mu_{(3,2)}\geq 0$, we have
$\sum_{i=1}^{3}\bigl(\alpha_{3}(x)-\alpha_{i}(x)\bigr)
\varphi_{i}(x)\mu_{(3,2)}\leq 0$, and on the other hand
$\eta_{2}\varphi_{2}(x)\sum_{i=1}^{3}\varphi_{i}(x)<0$ so that,
if $\sum_{i=1}^{3}\bigl(\alpha_{2}(x)-\alpha_{i}(x)\bigr) \varphi_{i}(x)$
happens to be nonnegative, then we also have
$\sum_{i=1}^{3}\bigl(\alpha_{2}(x)-\alpha_{i}(x)\bigr)
\varphi_{i}(x)\mu_{(2,2)}\leq 0$, and then we are done since
$\beta_{2}$, being a sum of negative terms (one of them strictly
negative) must be negative as well. Otherwise, $$
\sum_{i=1}^{3}\Bigl(\alpha_{2}(x)-\alpha_{i}(x)\Bigr)
\varphi_{i}(x) $$ is negative, and then we finally pass to the
term $\beta_{1}$ multiplying $(x-y_{1})$ in $(14)$, that is,
    $$
    \beta_{1}:=\eta_{1}\varphi_{1}(x)\sum_{i=1}^{3}\varphi_{i}(x) +
    \sum_{k=1}^{3}\Biggl(\sum_{i=1}^{3}\Bigl(\alpha_{k}(x)-
    \alpha_{i}(x)\Bigr)\varphi_{i}(x)\Biggr)
    \mu_{(k,1)}.
    $$
Here, by the assumptions we have made so far and taking into
account the signs of $\mu_{(k,j)}$ and $\eta_{j}$, we see that
$\sum_{i=1}^{3}\bigl(\alpha_{k}(x)- \alpha_{i}(x)\bigr)\varphi_{i}(x)
\mu_{(k,1)}\leq 0$ for $k=2, 3$. Having arrived at this point, it
is sure that
$\sum_{i=1}^{3}\bigl(\alpha_{1}(x)-\alpha_{i}(x)\bigr)\varphi_{i}(x)$ must
be nonnegative (otherwise the numbers
$\sum_{i=1}^{3}\bigl(\alpha_{k}(x)-\alpha_{i}(x)\bigr)\varphi_{i}(x)$
should be strictly negative for all $k=1,2,3$, which is
impossible if one takes $\alpha_{k}(x)$ to be the maximum of the
$\alpha_{i}(x)$), and now we can deduce as before that
$\beta_{1}<0$.

Finally let us consider the case when some of the $\varphi_{i}(x)$
vanish, for $i=1,2$ (remember that $\varphi_{3}(x)\neq 0$ since
$x\in B_{3}$, the open support of $\varphi_{3}$). From the
definitions of $\mu_{(k,j)}$, $g_{n}$ and $\varphi_{n}$, it is
clear that $\mu_{(k,j)}=0$ whenever $\varphi_{j}(x)=0$ or
$\varphi_{k}(x)=0$, and bearing this fact in mind we can simplify
equality $(14)$ to a great extent by dropping all the terms that
now vanish.

If
$\varphi_{1}(x)=\varphi_{2}(x)=0$ then $(14)$ reads
    $$
    \varphi_{3}(x)^{2}\eta_{3}\,(x-y_{3})=0,
    $$
which cannot happen since we assumed $x\neq y_{j}$ (this means
that the only critical point that $f_{n}$ can have in
$B_{3}\setminus(B_{1}\cup B_{2})$ is $y_{3}$).

If
$\varphi_{1}(x)=0$ and $\varphi_{2}(x)\neq 0$ then the term
$\beta_{1}$ accompanying $(x-y_{1})$ in $(14)$ vanishes, and hence
$(14)$ is reduced to
    $$
    \sum_{j=2}^{3}\Biggl[
    \eta_{j}\varphi_{j}(x)\sum_{i=2}^{3}\varphi_{i}(x) +
    \sum_{k=j}^{3}\Bigl(\sum_{i=2}^{3}(\alpha_{k}(x)-
    \alpha_{i}(x))\varphi_{i}(x)\Bigr)
    \mu_{(k,j)}\Biggr]\,(x-y_{j})=0.
    $$
Since at least one of the numbers $\sum_{i=2}^{3}(\alpha_{k}(x)-
\alpha_{i}(x))\varphi_{i}(x)$, $k=2,3$, is nonnegative, the same
reasoning as in the first case allows us to conclude that either
$\beta_{3}$ or $\beta_{2}$ is strictly negative. Finally, in the
case $\varphi_{1}(x)\neq 0$ and $\varphi_{2}(x)=0$, it is
$\beta_{2}$ that vanishes, and $(14)$ reads $\beta_{1}\,(x-y_{1})+
\beta_{3}\,(x-y_{3})=0$, where
    $$
    \beta_{3}=\eta_{3}\varphi_{3}(x)\sum_{i=1, i\neq 2}^{3}\varphi_{i}(x) +
    \sum_{i=1, i\neq 2}^{3}\Bigl(\alpha_{3}(x)-
    \alpha_{i}(x)\Bigr)\varphi_{i}(x)
    \mu_{(3,3)},
    $$
and $$
    \beta_{1}=\eta_{1}\varphi_{1}(x)\sum_{i=1, i\neq 2}^{3}\varphi_{i}(x) +
    \sum_{k=1, i\neq 2}^{3}\sum_{i=1, i\neq 2}^{3}\Bigl(\alpha_{k}(x)-
    \alpha_{i}(x)\Bigr)\varphi_{i}(x)
    \mu_{(k,1)}.
    $$
Again, at least one of the numbers $\sum_{i=1, i\neq
2}^{3}\bigl(\alpha_{k}(x)- \alpha_{i}(x)\bigr)\varphi_{i}(x)$, $k=1,3$, is
nonnegative, and the same argument as above applies.

To finish the proof of the proposition we will need even more
accurate information about the location of the critical points of
$f_{n}$. Bearing in mind the definition of the functions
$\varphi_{j}$, whose open support are the $B_{j}$, it is clear
that the above discussion shows, in fact, the following
inclusions:
\begin{itemize}
\item[{}] $C_{f_{3}}\cap B_{3}\subseteq\mathcal{A}[y_{1}, y_{2},
y_{3}]$;
\item[{}] $C_{f_{3}}\cap(B_{3}\setminus B_{1})\subseteq\mathcal{A}[y_{2},
y_{3}]$, and $C_{f_{3}}\cap(B_{3}\setminus
B_{2})\subseteq\mathcal{A}[y_{1}, y_{3}]$ ;
\item[{}] $C_{f_{3}}\cap(B_{3}\setminus (B_{1}\cup
B_{2}))\subseteq \mathcal{A}[y_{3}]$.
\end{itemize}
An analogous argument in the case $n\geq 4$ proves the second part
of the statement of Fact~\ref{the critical set is locally finite
dimensional}.
\end{proof}

\begin{rem}
{\em Note that from Fact~\ref{the critical set is locally finite
dimensional} it follows that the set of critical points
$C_{\varphi}$ is locally compact, since it is closed and it is
locally a bounded set of a finite-dimensional affine subspace.}
\end{rem}

\medskip

So far, all the properties we have shown about our functions
$f_{n}$ are independent of the way we may choose the numbers
$\lambda_{j}$ in the definitions of $B_j$ and $\varphi_{j}$. Now
we are going to be more accurate and see how we can select those
numbers $\lambda_{j}$ so as to have more control over the set
$C_{\varphi}$ of critical points of $\varphi$. Indeed, we want
$C_{\varphi}$ not only to be locally compact, but to consist of a
sequence of suitably isolated small compact sets $K_{n}$. That
is, we want to write
    $
    C_{\varphi}\subseteq\bigcup_{n=1}^{\infty}K_{n},
    $
where the $K_n$ are compact sets which are associated with open
sets $U_{n}$ so that $K_{n}\subset U_{n}\subset B(y_{n}, r_{n})$,
and $U_{n}\cap U_{m}=\emptyset$ whenever $n\neq m$.

We will choose the numbers $\lambda_{n}$ and the open sets $U_{n}$
inductively.

\medskip
\noindent {\bf First step.} Define $\varphi_{1}$ as above and put
$f_{1}(x)=\alpha_{1}(x)$ for all $x\in B_{1}=B(y_{1}, r_{1})$.
Set $\mu_{2}=1/2$, $K_{1}=C_{f_{1}}\cap B_{1}=\{y_{1}\}$, and
$U_{1}=B(y_{1},\mu_{2}r_{1})$.

\medskip
\noindent {\bf Second step.} Fix $\lambda_{2}\in (\mu_{2}, 1)$,
and define $B_{2}$, $\varphi_{2}$, and $f_{2}$ as above.
According to Fact \ref{the critical set is locally finite
dimensional}, we have that
\begin{itemize}
\item[{}] $C_{f_{2}}\cap B_{2}\subset\mathcal{A}[y_{1}, y_{2}]$,
and
\item[{}] $C_{f_{2}}\cap (B_{2}\setminus B_{1})
\subseteq\mathcal{A}[y_{2}]$.
\end{itemize}
We claim that there must exist some $\mu_{3}\in (\lambda_{2}, 1)$
so that $\overline{C_{f_{2}}\cap B_{2}\cap B_{1}}\subset B(y_{1},
\mu_{3}r_{1})$. Otherwise there would exist a sequence $(x_{j})$
in $C_{f_{2}}\cap B_{2}\cap B_{1}$ so that $\|x_{j}-y_{1}\|$ goes
to $r_{1}$ as $j$ goes to $\infty$. Since $C_{f_{2}}\cap
B_{2}\subset\mathcal{A}[y_{1}, y_{2}]$, we may assume, by
compactness, that $x_{j}$ converges to some point
$x_{0}\in\partial B(y_{1}, r_{1})=S_{1}$. If $x_{0}\in B(y_{2},
r_{2})$ then $f_{2}'(x_{0})=0$ (by continuity of $f_{2}'$), and
$x_{0}\neq y_{2}$ (because $y_{2}\notin S_{1}$ by ii) of Lemma
\ref{controlled intersection of spheres}), so
    $$
    f_{2}'(x_{0})=\alpha_{2}'(x_{0})\neq 0,
    $$
a contradiction. Therefore it must be the case that
$x_{0}\in\partial B(y_{2}, r_{2})=S_{2}$. But then
    $$
    x_{0}\in S_{1}\cap S_{2}\cap\mathcal{A}[y_{1}, y_{2}],
    $$
and this contradicts Lemma \ref{controlled intersection of
spheres}.

So let us take $\mu_{3}\in (\lambda_{2}, 1)$ so that
$\overline{C_{f_{2}}\cap B_{2}\cap B_{1}}\subset B(y_{1},
\mu_{3}r_{1})$. In the case that $y_{2}\in B_{1}$, let us simply
set
\begin{itemize}
\item[{}] $U_{2}=B(y_{2}, r_{2})\cap B(y_{1}, \mu_{3}r_{1})\setminus
\overline{B}(y_{1}, \mu_{2}r_{1})$, and
\item[{}] $K_{2}=\overline{C_{f_{2}}\cap B_{2}\cap B_{1}}\subset
U_{2}$.
\end{itemize}
In the case that $y_{2}\notin B_{1}$, find $\delta_{2}\in (0,
\mu_{3} r_{2})$ so that $B(y_{2}, \delta_{2})\subset
B_{2}\setminus\overline{B_{1}}$, and set
\begin{itemize}
\item[{}] $U_{2}=\big[B(y_{2}, r_{2})\cap B(y_{1}, \mu_{3}r_{1})\setminus
\overline{B}(y_{1}, \mu_{2}r_{1})\big]\cup B(y_{2},\delta_{2})$,
and
\item[{}] $K_{2}=\overline{C_{f_{2}}\cap B_{2}\cap B_{1}}\cup\{y_{2}\}\subset
U_{2}$.
\end{itemize}
Clearly, we have that $C_{f_{2}}\subseteq K_{1}\cup K_{2}$, and
$U_{1}\cap U_{2}=\emptyset$.

\medskip
\noindent {\bf Third step.} Now choose $\lambda_{3}\in (\mu_{3},
1)$ with $\lambda_{3}>1-1/3$, and define $B_{3}$, $\varphi_{3}$,
and $f_{3}$ as above. We have that $f_{3}$ and $f_{2}$ coincide on
$(B_{1}\cup B_{2})\setminus B_{3}$. On $B_{3}$, according to Fact
\ref{the critical set is locally finite dimensional}, we know that
\begin{itemize}
\item[{}] $C_{f_{3}}\cap B_{3}\cap B_{2}\cap B_{1}\subseteq\mathcal{A}[y_{1}, y_{2},
y_{3}]$;
\item[{}] $C_{f_{3}}\cap(B_{3}\cap B_{2}\setminus B_{1})\subseteq\mathcal{A}[y_{2},
y_{3}]$, and $C_{f_{3}}\cap(B_{3}\cap B_{1}\setminus
B_{2})\subseteq\mathcal{A}[y_{1}, y_{3}]$ ;\hfill $(15)$
\item[{}] $C_{f_{3}}\cap(B_{3}\setminus (B_{1}\cup
B_{2}))\subseteq \mathcal{A}[y_{3}]$.
\end{itemize}
Again, there must be some $\mu_{4}\in (\lambda_{3}, 1)$ so that
    $$
    \overline{C_{f_{3}}\cap B_{3}\cap(B_{1}\cup B_{2})}
    \subset B(y_{1}, \mu_{4}r_{1})\cup
    B(y_{2}, \mu_{4}r_{2}).
    $$
Otherwise (bearing in mind the local compactness of
$\mathcal{A}[y_{1}, y_{2}, y_{3}]$), there would exist a sequence
$(x_{j})$ in $C_{f_{3}}\cap B_{3}\cap (B_{1}\cup B_{2})$ so that
$(x_{j})$ converges to some point $x_{0}$ and $(x_{j})$ is not
contained in $B(y_{1}, \mu_{4}r_{1})\cup
    B(y_{2}, \mu_{4}r_{2})$ for any $\mu_{4}<1$. Since a subsequence
of $(x_{j})$ must be contained in one of the sets listed in
$(15)$, we deduce that the limit point $x_{0}$ must belong to one
of the following sets:
\begin{itemize}
\item[{}] $S_{2}\cap S_{1}\cap\mathcal{A}[y_{1}, y_{2},
y_{3}]$;
\item[{}] $S_{2}\cap\mathcal{A}[y_{2},
y_{3}]\setminus B_{1}$;
\item[{}] $S_{1}\cap\mathcal{A}[y_{1},
y_{3}]\setminus B_{2}$,
\end{itemize}
Now we have two cases: either $x_{0}\in B_{3}$, or $x\in\partial
B_{3}$. If $x_{0}\in B_{3}$ then $f_{3}'(x_{0})=0$ (by continuity
of $f_{3}'$), and $x_{0}\neq y_{3}$ (because $y_{3}\notin
S_{1}\cup S_{2}$ by (ii) of Lemma~\ref{controlled intersection of
spheres}), so it follows that
    $$
    f_{3}'(x_{0})=\alpha_{3}'(x_{0})\neq 0,
    $$
a contradiction. On the other hand, if $x_{0}\in\partial B_{3}$
then $x_{0}\in S_{3}$ as well, and now one of the following must
hold:
\begin{itemize}
\item[{}] $x_{0}\in S_{3}\cap S_{2}\cap S_{1}\cap\mathcal{A}[y_{1}, y_{2},
y_{3}]$;
\item[{}] $x_{0}\in S_{3}\cap S_{2}\cap\mathcal{A}[y_{2},
y_{3}]$;
\item[{}] $x_{0}\in S_{3}\cap S_{1}\cap\mathcal{A}[y_{1},
y_{3}]$,
\end{itemize}
but in any case this contradicts Lemma \ref{controlled
intersection of spheres}.

Hence we can take $\mu_{4}\in (\lambda_{3}, 1)$ so that
    $$
    \overline{C_{f_{3}}\cap B_{3}\cap(B_{1}\cup B_{2})}
    \subset B(y_{1}, \mu_{4}r_{1})\cup
    B(y_{2}, \mu_{4}r_{2}).
    $$
Now two possibilities arise. If $y_{3}\in B_{1}\cup B_{2}$, let
us define
    $$
    U_{3}=\Biggl[B(y_{3}, r_{3})\setminus\bigcup_{j=1}^{2}
    \overline{B}(y_{j}, \mu_{3}r_{j})\Biggr]\bigcap
    \Biggl[\bigcup_{j=1}^{2}B(y_{j}, \mu_{4}r_{j})\Biggr],
    $$
and
    $$
    K_{3}=\overline{C_{f_{3}}\cap B_{3}\cap(B_{1}\cup B_{2})}\subset U_{3}.
    $$
If $y_{3}\notin B_{1}\cup B_{2}$, since $y_{3}\notin S_{1}\cup
S_{2}$ we can find $\delta_{3}\in (0, \mu_{4}r_{3})$ so that
$B(y_{3},\delta_{3})\subseteq B_{3}\setminus (B_{1}\cup B_{2})$,
and then we can set
    $$
    U_{3}=\Biggl[\Bigl(B(y_{3}, r_{3})\setminus\bigcup_{j=1}^{2}
    \overline{B}(y_{j}, \mu_{3}r_{j})\Bigr)\bigcap
    \Bigl(\bigcup_{j=1}^{2}B(y_{j}, \mu_{4}r_{j})\Bigr)\Biggr]\bigcup
    B(y_{3}, \delta_{3}),
    $$
and
    $$
    K_{3}=\overline{[C_{f_{3}}\cap B_{3}\cap(B_{1}\cup B_{2})]\cup
    \{y_{3}\}}\subset U_{3}.
    $$
Notice that $U_{3}$ does not meet $U_{1}$ or $U_{2}$, and
$C_{f_{3}}\subseteq K_{1}\cup K_{2}\cup K_{3}$.

\medskip
\noindent {\bf N-th step.} Suppose now that $\mu_{j}$,
$\lambda_{j}$, $\varphi_{j}$, $B_{j}$, $f_{j}$, $K_{j}$, $U_{j}$ have already
been fixed for $j=1, ..., n$ (and also $\mu_{n+1}$ has been
chosen) in such a manner that $f_{j}$ agrees with $f_{j-1}$ on
$(B_{1}\cup ...\cup B_{j-1})\setminus B_{j}$, and $K_{j}$ and
$U_{j}$ are of the form
    $$
    K_{j}=\overline{C_{f_{j}}\cap B_{j}\cap(B_{1}\cup... \cup B_{j-1})}
    \eqno(16)
    $$
and
    $$
    U_{j}=\biggl[B(y_{j}, r_{j})\setminus\Bigl(\bigcup_{i=1}^{j-1}
    \overline{B}(y_{i}, \mu_{j}r_{i})\Bigr)\biggr]\bigcap
    \Bigl[\bigcup_{i=1}^{j-1}B(y_{i}, \mu_{j+1}r_{i})\Bigr] \eqno(17)
    $$
in the case that $y_{j}\in B_{1}\cup ...\cup B_{j-1}$, and are of
this form plus $\{y_{j}\}$ and $B(y_{j},\delta_{j})$ respectively
when $y_{j}\notin B_{1}\cup ...\cup B_{j-1}$; assume additionally
that $U_{j}\cap U_{k}=\emptyset$ whenever $j\neq k$, that
$C_{f_{j}}\subseteq \bigcup_{i=1}^{j}K_{i}$, and that
$\lambda_{j}>1-1/j$. Let us see how we can choose $\lambda_{n+1}$,
$\mu_{n+2}$, $K_{n+1}$ and $U_{n+1}$ so that the extended bunch
keeps the required properties.

Pick any $\lambda_{n+1}\in (\mu_{n+1}, 1)$ so that
$\lambda_{n+1}>1-1/(n+1)$, and define $\varphi_{n+1}$, $B_{n+1}$
and $f_{n+1}$ as above. We know that $f_{n+1}$ agrees with
$f_{n}$ on the set $(B_{1}\cup ... \cup B_{n})\setminus B_{n+1}$.
On $B_{n+1}$, according to Fact \ref{the critical set is locally
finite dimensional}, we have that
    $$
    C_{f_{n+1}}\cap\big(B_{n+1}\setminus\bigcup_{j=1}^{m}B_{k_{j}}\big)
    \subseteq\mathcal{A}\big[\{y_{1}, ..., y_{n+1}\}
    \setminus\{y_{k_{1}}, ..., y_{k_{m}}\}\big]
    $$
for every finite sequence of integers $0<k_{1}<k_{2}< ...
<k_{m}<n+1$.

We claim that there exists some $\mu_{n+2}\in (\lambda_{n+1}, 1)$
so that
    $$
    \overline{C_{f_{n+1}}\cap B_{n+1}\cap (B_{1}\cup ... \cup B_{n})}
    \subseteq \bigcup_{i=1}^{n}B(y_{i},\mu_{n+2}r_{i}).
    $$
Otherwise there would exist a finite (possibly empty!) sequence of
integers $0<k_{1}<k_{2}< ... <k_{m}<n+1$, and a sequence
$(x_{j})_{j=1}^{\infty}$ contained in
    $$
    \Biggl[C_{f_{n+1}}\cap B_{n+1}\cap \Bigl(\bigcap_{j=1}^{\ell}B_{i_{j}}\Bigr)\Biggr]
    \setminus\Bigl(\bigcup_{j=1}^{m}B_{k_{j}}\Bigr)
    \subseteq\mathcal{A}[y_{i_{1}}, ..., y_{i_{\ell}}, y_{n+1}]
    $$
(where $i_{1}, ...,i_{\ell}$ are the positive integers less than
or equal to $n$ that are left when we remove $k_{1}, ...,
k_{m}$), such that $(x_{j})$ converges to some point $x_{0}\in
S_{i_{1}}\cap ...\cap S_{i_{\ell}}$ with
$x_{0}\notin\bigcup_{j=1}^{m}B_{k_{j}}$.

If $x_{0}\in B_{n+1}$ then $f_{n+1}'(x_{0})=0$ (by continuity of
$f_{n+1}'$), and $x_{0}\neq y_{n+1}$, so we easily see that
    $$
    f_{n+1}'(x_{0})=\alpha_{n+1}'(x_{0})\neq 0,
    $$
a contradiction.

If $x_{0}\in \partial B_{n+1}$ then $x_{0}\in S_{n+1}$ as well,
and in this case we have
    $$
    x_{0}\in S_{i_{1}}\cap ...\cap S_{i_{\ell}}\cap S_{n+1}\cap
    \mathcal{A}[y_{i_{1}}, ..., y_{i_{\ell}}, y_{n+1}],
    $$
but this contradicts Lemma \ref{controlled intersection of
spheres}.

Therefore we may take $\mu_{n+2}\in (\lambda_{n+1}, 1)$ so that
    $$
    \overline{C_{f_{n+1}}\cap B_{n+1}\cap (B_{1}\cup... \cup
    B_{n})}\subseteq
    \bigcup_{i=1}^{n}B(y_{i},\mu_{n+2}r_{i}).
    $$
As before, now we face two possibilities. If
$y_{n+1}\in\bigcup_{i=1}^{n}B_{i}$, let us define
    $$
    U_{n+1}=\Biggl[B(y_{n+1}, r_{n+1})\setminus\bigcup_{i=1}^{n}
    \overline{B}(y_{i}, \mu_{n+1}r_{i})\Biggr]\bigcap
    \Biggl[\bigcup_{i=1}^{n}B(y_{i}, \mu_{n+2}r_{i})\Biggr],
    $$
and $$
    K_{n+1}=\overline{C_{f_{n+1}}\cap B_{n+1}\cap(B_{1}\cup... \cup
    B_{n})}.
    $$
If $y_{n+1}\notin\bigcup_{i=1}^{n}B_{i}$, since $y_{n+1}\notin
S_{i}$ we may find $\delta_{n+1}\in (0, \mu_{n+2}r_{n+1})$ so
that $B(y_{n+1},\delta_{n+1})\subseteq B_{n+1}\setminus
\bigcup_{i=1}^{n}B_{i}$, and then we can add this ball to the
above $U_{n+1}$, and the point $\{y_{n+1}\}$ to that $K_{n+1}$,
in order to obtain sets $U_{n+1}$, $K_{n+1}$ with the required
properties.

By induction, the sequences $(\varphi_{n})$, $(f_{n})$, $(U_{n})$,
$(K_{n})$ are well defined and satisfy the above properties.

From the construction it is clear that $U_{n}\cap U_{m}=\emptyset$
whenever $n\neq m$, and
    $$
    C_{f_{n}}\subseteq\bigcup_{j=1}^{n}K_{j}
    $$
for all $n$. Note also that $U_{n}\subseteq B(y_{n}, r_{n})$ for
all $n$, and $\lim_{n\to\infty}\lambda_{n}=1$.

As observed before,
    $$
    \varphi(x)=\frac{\sum_{n=1}^{\infty}\alpha_{n}(x)\varphi_{n}(x)}
    {\sum_{n=1}^{\infty}\varphi_{n}(x)}=\lim_{n\to\infty}f_{n}(x)
    $$
and, moreover, for each $x\in X$ there exists an open
neighborhood $V_{x}$ of $x$ and some $n_{x}\in\mathbb{N}$ so that
$\varphi(y)=f_{n_{x}}(y)$ for all $y\in V_{x}$. Bearing these
facts in mind, it is immediately checked that
$C_{\varphi}\subseteq\bigcup_{n=1}^{\infty}K_{n}$. Now it is
clear that $\varphi$ satisfies (a), (b) and (d) in the statement
of Proposition \ref{existence of varphi}. On the other hand,
remember that (c) is a consequence of fact \ref{good
approximation}.


\medskip
\begin{rem}\label{remark for e(x)}
{\em Let us say a few words as to the way one has to modify the
above proofs in order to establish Theorem \ref{main theorem} when
$\varepsilon$ is a positive continuous function. At the beginning
of the proof of Proposition \ref{existence of varphi}, before
choosing the $\delta_{x}$, we have to take some number
$\alpha_{x}>0$ so that
$|\varepsilon(y)-\varepsilon(x)|\leq\varepsilon(x)/2$ whenever
$\|y-x\|\leq 2\alpha_{x}$ and then we can find some
$\delta_{x}\leq\alpha_{x}$ so that
$|f(y)-f(x)|\leq\varepsilon(x)/4$ whenever $y\in B(x,
2\delta_{x})$. Equation $(7)$ above reads now
    $$
    |f(y)-f(y_{n})|\leq\varepsilon(y_{n})/2
    $$
for all $y\in B(y_{n}, r_{n})$. Some obvious changes must be made
in the definition of the functions $a_{n}$ and $\alpha_{n}$. Fact
\ref{good approximation} and Proposition \ref{existence of
varphi}(c) can be reduced to saying that $$|\varphi(y)-f(x)|\leq
2\varepsilon(y_{n})$$ for all $x, y\in B(y_{n}, r_{n})$ and each
$n\in\mathbb{N}$. Finally, at the end of the proof of Theorem
\ref{main theorem} we get that
    $$
    |\psi(x)-f(x)|\leq 2\varepsilon(y_{n})
    $$
whenever $x, y\in B(y_{n}, r_{n})$; now, taking into account that
$r_{n}\leq\alpha_{y_{n}}$, we have that $\varepsilon(y_{n})\leq
2\varepsilon(x)$ for all $x\in B(y_{n}, r_{n})$. Hence, by
combining these inequalities, we obtain that
    $
    |\psi(x)-f(x)|\leq 4\varepsilon(x)
    $
for all $x\in X$.}
\end{rem}



\medskip
\begin{center}
{\bf Proof of Theorem \ref{removing compact sets}}
\end{center}

The proof of Theorem~\ref{removing compact sets} is done in two
steps. The first one uses the noncomplete norm technique of
deleting compact sets introduced in \cite{ADo, Do}. We only
sketch the guidelines of this part, referring to the proof of
Theorem~2.1 in \cite{ADo} for the details. We will show that a
mapping of the form $G(x)=x+p(f(x))$, $x\in X\setminus K$, for a
certain function $f:X\to[0,+\infty)$ with $f^{-1}(K)=0$ and a path
$p:(0,+\infty)\to X$, establishes a $C^\infty$ diffeomorphism
between $X\setminus K$ and $X$. The map $G$ can be viewed as a
{\em small} perturbation of the identity. In order that the
perturbation $p\circ f$ be small, $p$ and $f$ must satisfy some
Lipschitzian-type conditions with respect to a certain distance
induced by a smooth {\em noncomplete norm} $\omega$. Lemma
\ref{Whitney function for compacta} provides us with a required
function $f(x)$ which can be viewed as a smooth substitute for
the $\omega$-distance function from $x$ to the set $K$. Lemma
\ref{deleting path} gives us a required path $p(t)$ which avoids
compact sets and gets lost in the infinitely many dimensions of
$X$ as $t$ goes to $0$; by pushing away $\omega$-neighborhoods
of $K$ along the path $p$, the mapping $G^{-1}$ will make $K$
disappear. By combining all these tools, the $C^\infty$
diffeomorphism $G$ can be constructed in such a way that $G$
restricts to the identity outside a given $\omega$-neighborhood
of $K$.

So far this is the same negligibility scheme as in \cite{ADo,
Do}. The second step of the proof is to construct a
self-diffeomorphism $F$ of $X$ that fixes the compact set $K$ and
takes the open set $U$ (which in general is {\em not} a
$\omega$-neighborhood of $K$) onto a $\omega$-neighborhood of
$K$, and then to adjust the definition of $G$ so that it
restricts to the identity outside $F(U)$. If we succeed in doing
so then the composition $h=F^{-1}\circ G^{-1}\circ F$ will define
a diffeomorphism from $X$ onto $X\setminus K$ with the property
that $h$ is the identity outside $U$.

\begin{defn}\label{omega neighbourhoods}
{\em Let $(X, \|\cdot\|)$ be a Banach space. We say that a norm
$\omega:X\longrightarrow[0,+\infty)$ is a $C^p$ smooth noncomplete
norm on $X$ provided $\omega$ is $C^p$ smooth (with respect to
$\|\cdot\|$) away from the origin, but the norm $\omega$ is not
equivalent to $\|\cdot\|$. Geometrically speaking, this means
that the unit ball of $\omega$ is a symmetric $C^p$ smooth convex
body that contains no rays and yet is unbounded. We define the
(open) $\omega$-ball of center $x$ and radius $r$ as
    $$
    B_{\omega}(x,r)=\{y\in X : \omega(y-x)<r\},
    $$
and the $\omega$-distance from $x$ to $A$ as
    $$
    d_{\omega}(x,A)=\inf\{\omega(x-z) : z\in A\}.
    $$
We say that a set $V$ is an $\omega$-neighborhood of a subset $A$
of $X$ provided that for every $x\in A$ there exists some $r>0$
so that $B_{\omega}(x,r)\subseteq V$.}
\end{defn}

We next state the two facts we need for the first part of the
proof. All the omitted proofs can be found in \cite{ADo, Do}.

\begin{lem}\label{Whitney function for compacta}
Let $\omega:X\longrightarrow [0,+\infty)$ be a $C^\infty$ smooth
noncomplete norm in the Hilbert space $X$, and let $K$ be a
compact subset of $X$. Then, for each $\varepsilon>0$ there
exists a continuous function $f=f_\varepsilon:X\longrightarrow
[0,+\infty)$ such that
\begin{enumerate}
\item $f$ is $C^{\infty}$ smooth on $X\setminus K$;
\item $f(x)-f(y)\leq\omega(x-y)$ for every $x,
y\in X$;
\item $f^{-1}(0)=K$;
\item $\inf\{f(x)\mid d_{\omega}(x,K)\geq\eta\}>0$
for every $\eta>0$;
\item $f$ is constant on the set
$\{x\in X\mid d_{\omega}(x,K)\geq\varepsilon\}$.
\end{enumerate}
\end{lem}

\begin{lem}\label{deleting path}
Let $\omega$ be a continuous noncomplete norm in the Hilbert
space $X$. Then, for every $\delta>0$, there exists a
$C^{\infty}$ path $p=p_{\delta}:(0,+\infty)\longrightarrow X$ such
that
\begin{enumerate}
\item $\omega(p(\alpha)-p(\beta))\leq
\frac{1}{2}(\beta-\alpha)$ if $\beta\ge\alpha>0$;
\item For every compact set $A\subset
X$ there exists $t_{0}>0$ such that $$ \inf\{\omega(z-p(t))\mid
0<t\leq t_{0}, z\in A\} >0; $$
\item $p(t)=0$ if and only if $t\geq\delta$.
\end{enumerate}
\end{lem}
The following lemma is the key to the second step of the proof,
allowing us to improve, at least for the Hilbert case, the
negligibility scheme introduced in \cite{ADo, Do}. Note also that
the norm $\omega$ that we will use in the first step is in fact
the one provided by this lemma.
\begin{lem}\label{radial push}
Let $(X, \|\cdot\|)$ be an infinite-dimensional Hilbert space
(with its usual hilbertian norm). Then, for every compact set $K$
and every open set $U$ so that $K\subset U$, there exist a
$C^\infty$ diffeomorphism $F:X\longrightarrow X$ and a $C^\infty$
smooth noncomplete norm $\omega$ on $X$ such that $F(K)=K$ and
$F(U)$ is an $\omega$-neighborhood of $K$.
\end{lem}
\begin{proof}
Since $K$ is compact and $U$ is an open neighborhood of $K$ we
can find points $x_{1}, ..., x_{n}\in K$ and positive numbers
$r_{1}, ..., r_{n}$ so that
    $$
    K\subset \bigcup_{i=1}^{n} B_{\|.\|}(x_{i}, r_{i}) \subset
    \bigcup_{i=1}^{n} B_{\|.\|}(x_{i}, 2r_{i})\subseteq U. \eqno
    (18)
    $$
We may assume that $0\in K$. Let $Y=\textrm{span}\{x_{1}, ...,
x_{n}\}$, and write $X=Y\oplus Z$, where $Z$ is an
infinite-dimensional space of finite codimension that is
orthogonal to $Y$. Since the norm $\|\cdot\|$ is hilbertian we
have that
    $$
    \|x\|=\|(y,z)\|=\big( \|y\|^{2}+\|z\|^{2}\big)^{1/2}
    $$
for every $x=(y,z)\in X=Y\oplus Z$.

Take a normalized basic sequence $(z_{i})$ in $Z$ so that the
vectors $z_{i}$ are pairwise orthogonal, and let $W$ be the closed
linear subspace spanned by $(z_i)$. Let us write $Z=W\oplus V$,
where $V$ is the orthogonal complement of $W$ in $Z$. Define
$\omega_{Z}:Z\longrightarrow [0,+\infty)$ by
    $$
    \omega_{Z}(w,v)=\biggl[\sum_{j=1}^{\infty}\Bigl(\frac{<w,z_{j}>}{2^{j}}\Bigr)^{2}
    +\|v\|^{2}\biggr]^{1/2},
    $$
where $<,>$ denotes the inner product on $X$. Then $\omega_{Z}$ is
a $C^\infty$ smooth noncomplete norm on $Z$, as it is easily
checked. We also have that $\omega_{Z}(z)\leq \|z\|$ for every
$z\in Z$. If we define now $\omega:X=Y\oplus Z\longrightarrow
[0,+\infty)$ by
    $$
    \omega(x)=\omega(y,z)=\big(\|y\|^{2}+\omega_{Z}(z)^{2}\big)^{1/2},
    $$
it is clear that $\omega$ is a $C^\infty$ smooth noncomplete norm
on $X$ (note that in fact both $\omega_{Z}$ and $\omega$ are
real-analytic, as they are prehilbertian).

For each $i=1, ..., n$, let us now pick $C^\infty$ smooth
functions $\theta_{i}:\mathbb{R}\longrightarrow [0,1]$ so that
$\theta_i$ is nondecreasing and $\theta_{i}^{-1}(0)=(-\infty,
r_{i}]$, while $\theta_{i}^{-1}(1)=[2r_{i}, +\infty)$. Define then
$g:X=Y\oplus Z \longrightarrow [0, 1]$ by
    $$
    g(y,z)=g(x)=\prod_{i=1}^{n}\theta_{i}(\|x-x_{i}\|)
    $$
for all $x\in X$. Note that the function $g$ is $C^\infty$ smooth
on $X$ and has the following properties:
\begin{itemize}
\item[{(i)}] the function $t\mapsto g(y, tz)$, $t\geq 0$, is nondecreasing,
for all $(y,z)\in X=Y\oplus Z$;
\item[{(ii)}] $g(x)=0$ if $x\in\bigcup_{i=1}^{n} B_{\|.\|}(x_{i}, r_{i})$;
\item[{(iii)}] $g(x)=1$ whenever $x\notin\bigcup_{i=1}^{n} B_{\|.\|}(x_{i},
2r_{i})$.
\end{itemize}
The first property is merely a consequence of the definition of
$g$ and the fact that the function $t\mapsto \|((y-x_{i}, tz)\|$,
$t\geq 0$, is increasing for every $(y,z)\in Y\oplus Z$ and every
$i=1, ..., n$. Note that here we are using that $\|\cdot\|$ is a
hilbertian norm; this property is not necessarily true for other
norms.

Let us define our mapping $F:X=Y\oplus Z\longrightarrow X$ by
    $$
    F(x)=F(y,z)=\biggl( y, \, \Bigl(g(x)\frac{\|z\|}{\omega_{Z}(z)}
    +1-g(x)\Bigr)\, z \, \biggr).
    $$
Clearly, $F$ is $C^\infty$ smooth. By using the facts that the
functions $t\mapsto g(y, tz)$, $t\geq 0$, are nondecreasing, and
that $\omega_{Z}(z)\leq \|z\|$ for every $z\in Z$, $y\in Y$, it is
not difficult to see that $F$ is a bijection from every ray $\{(y,
tz) : t\geq 0\}$ onto itself, and therefore $F$ is one-to-one
from $X$ onto $X$. Moreover, a standard application of the
implicit function theorem allows to show that $F^{-1}$ is
$C^\infty$ smooth as well, and hence $F$ is a diffeomorphism.

Finally, by the definitions of $g$ and $F$, it is clear that
    $
    F(K)=K
    $.
In fact, $F$ restricts to the identity on the set $\bigcup_{i=1}^{n}
B_{\|.\|}(x_{i}, r_{i})$, which contains $K$, because $g$ takes
the value $0$ on this set.

On the other hand, if
$x\notin\bigcup_{i=1}^{n} B_{\|.\|}(x_{i}, 2r_{i})$ we have
$g(x)=1$, so $F(y,z)=(y, \frac{\|z\|}{\omega(z)} z)$, and
therefore
\begin{eqnarray*}
  & &\omega(F(x)-x_{j})=\omega\Bigl(y-x_{j}, \frac{\|z\|}{\omega(z)}z\Bigr)=
  \biggl(\|y-x_{j}\|^{2}+\omega\Bigl(\frac{\|z\|}{\omega(z)}z\Bigr)^{2}\biggr)^{1/2}\\
  & &=\bigl(\|y-x_{j}\|^{2}+\|z\|^{2}\bigr)^{1/2}=\|x-x_{j}\|\geq 2r_{j}
\end{eqnarray*}
for each $j=1, ..., n$, which means that $F(x)\notin\cup_{i=1}^{n}
B_{\omega}(x_{i}, 2r_{i})$. Therefore, considering $(18)$, and
bearing in mind that, since $\omega(x)\leq\|x\|$, the
$\omega$-balls are larger than the $\|\cdot\|$-balls, we deduce
that
    $$
    K\subset \bigcup_{i=1}^{n} B_{\omega}(x_{i}, r_{i})
    \subset \bigcup_{i=1}^{n} B_{\omega}(x_{i}, 2r_{i})\subseteq
    F\big(\bigcup_{i=1}^{n} B_{\|.\|}(x_{i}, 2r_{i})\big)\subseteq F(U);
    $$
in particular we see that $F(U)$ includes a finite union of
$\omega$-balls which in turn includes K, and this shows that
$F(U)$ is a $\omega$-neighborhood of $K$.
\end{proof}

Let us now see how we can finish the proof of Theorem
\ref{removing compact sets}. First, for the given sets $K\subset
U$, take a non-complete norm $\omega$ and a diffeomorphism
$F:X\longrightarrow X$ with the properties of Lemma \ref{radial
push}. Since $F(U)$ is a $\omega$-neighborhood of $K$ and $K$ is
also compact in $(X,\omega)$, we can write
    $$
    K\subset \bigcup_{i=1}^{n}B_{\omega}(x_{i}, r_{i})\subseteq
    \bigcup_{i=1}^{n}B_{\omega}(x_{i}, 2r_{i})\subseteq F(U)
    $$
for some points $x_{1}, ..., x_{n}\in K$ and positive numbers
$r_{1}, ..., r_{n}$ (in fact such an expression appears in the
proof of \ref{radial push}). This in turn implies that
$d_{\omega}(x, K)\geq \min\{r_{1}, ..., r_{n}\}>0$ whenever $x\in
X\setminus F(U)$, as it is easily seen.

Now, for $\varepsilon=\min\{r_{1}, ..., r_{n}\}$, we can choose a
function $f=f_{\varepsilon}$ satisfying the properties of Lemma
\ref{Whitney function for compacta} (for the already selected
$\omega$). Assuming $f(x)=\delta>0$ whenever
$d_{\omega}(x,K)\geq\varepsilon$, select a path $p=p_{\delta}$
from Lemma \ref{deleting path}. With these choices, for every
$x\in X\setminus K$, define
    $$
    G(x)=x+p(f(x)).
    $$
Exactly as in the proof of Theorem 2.1 in \cite{ADo}, it can be
checked that $G$ is a $C^\infty$ diffeomorphism from $X\setminus
K$ onto $X$, with the property that $G(x)=x$ whenever
$d_{\omega}(x, K)\geq\varepsilon$. In particular, since
$d_{\omega}(x, K)\geq \varepsilon=\min\{r_{1}, ..., r_{n}\}$
whenever $x\in X\setminus F(U)$, we have that $G$ restricts to
the identity outside $F(U)$.

Finally, let us define $h=F^{-1}\circ G^{-1}\circ F$. Taking into
account the properties of the diffeomorphisms $F:X\longrightarrow
X$ and $G:X\setminus K\longrightarrow X$, it is clear that $h$ is
a $C^\infty$ diffeomorphism from $X$ onto $X\setminus K$ so that
$h$ is the identity outside $U$. 

\medskip

\begin{center}
{\bf Acknowledgements}
\end{center}
\noindent We wish to thank Pilar Cembranos, Jos\'e Mendoza,
Tijani Pakhrou and Ra\'ul Romero, who helped us to realize that
Fact \ref{the critical set is locally finite dimensional} fails
when the norm is not hilbertian.


\vspace{3mm} \noindent Departamento de An\'alisis Matem\'atico.
Facultad de Ciencias Matem\'aticas. Universidad Complutense.
28040 Madrid, SPAIN\\ \noindent Departamento de An\'alisis
Matem\'atico. Universidad de Sevilla. Sevilla, SPAIN. \noindent
{\em E-mail addresses:} daniel\_azagra@mat.ucm.es, mcb@us.es
\end{document}